\documentclass[10pt]{amsart}

\usepackage{amsfonts}
\usepackage{amscd}
\usepackage{amssymb}
\usepackage{amsthm}
\usepackage{amsmath}
\usepackage{graphicx}
\usepackage{color}
\usepackage{hyperref}

\newtheorem{thm}{Theorem}[section]
\newtheorem{lem}[thm]{Lemma}
\newtheorem{cor}[thm]{Corollary}
\newtheorem{prop}[thm]{Proposition}

\theoremstyle{definition}
\newtheorem{defn}[thm]{Definition}
\newtheorem{rem}[thm]{Remark}

\newcommand{\commentout}[1]{}

\def\bu{_\bullet}

\def\ot{\leftarrow}

\def\xrarrow{\xrightarrow} 
\def\xlarrow{\xleftarrow} 

\def\noteq{\neq}

\def\<{\left<}
\def\>{\right>}

\DeclareMathOperator{\End}{End}

\DeclareMathOperator{\simp}{simp}

\DeclareMathOperator{\Hom}{Hom}%
\DeclareMathOperator{\HOM}{HOM}

\newcommand{\field}[1]{\mathbb{#1}}
\newcommand{\ZZ}{\ensuremath{{\field{Z}}}}
\newcommand{\CC}{\ensuremath{{\field{C}}}}
\newcommand{\RR}{\ensuremath{{\field{R}}}}
\newcommand{\QQ}{\ensuremath{{\field{Q}}}}

\def\LL{\Lambda}

\newcommand{\cC}{\ensuremath{{\mathcal{C}}}}

\newcommand{\cG}{\ensuremath{{\mathcal{G}}}}

\newcommand{\cN}{\ensuremath{{\mathcal{N}}}}

\newcommand{\cU}{\ensuremath{{\mathcal{U}}}}
\newcommand{\cV}{\ensuremath{{\mathcal{V}}}}

\newcommand{\Wh}{\ensuremath{{\mathcal{W}\textit{h}}}}
\newcommand{\cX}{\ensuremath{{\mathcal{X}}}}

\def\a{\alpha}

\def\g{\gamma}
\def\d{\partial}

\def\f{\varphi}

\def\r{\rho}
\def\s{\sigma}
\def\Sig{\Sigma}
\def\t{\tau}

\def\z{\zeta}
\def\w{\omega}

\def\st{\,|\,}
\def\op{^{op}}

\begin{document}

\def\version{091129}

\title{Twisting Cochains and higher torsion}

\author{Kiyoshi Igusa}

\begin{abstract} 
This paper gives a short summary of the central role played by Ed
Brown's ``twisting cochains'' in higher Franz-Reidemeister (FR)
torsion and higher analytic torsion. Briefly, any fiber bundle
gives a twisting cochain which is unique up to fiberwise homotopy
equivalence. However, when they are based, the difference between
two of them is a higher algebraic K-theory class measured by
higher FR torsion. Flat superconnections are also equivalent to
twisting cochains.
\end{abstract}

\dedicatory{This paper is dedicated to Edgar Brown.}

\address{Department of Mathematics, Brandeis University, Waltham, MA 02454}

\email{igusa@brandeis.edu}




\thanks{Research for this paper was supported by NSF}

\maketitle

\tableofcontents

%
%

\section*{Introduction}\label{sec:intro}

About 50 years ago Ed Brown \cite{Brown59:TwistedTensor}
constructed a small chain complex giving the homology of the total
space $E$ of a fiber bundle
\[
    F\to E\to B
\]
whose base $B$ and fiber $F$ are finite cell complexes. It is
given by the tensor product of chain complexes for $F$ and $B$
with the usual tensor product boundary map modified by a
``twisting cochain.'' There are many ways to understand the
meaning of the twisting cochain.
\begin{enumerate}
    \item It is the difference between two $A_\infty$ functors.
    \item It is a combinatorial flat $\ZZ$-graded superconnection.
    \item It is a family of chain complexes homotopy equivalent to $F$ and parametrized by $B$.
\end{enumerate}

If $F\to E\to B$ is a smooth bundle with compact manifold fiber
and simply connected base then we get another twisting cochain
given by fiberwise Morse theory. Comparison of these two twisting
cochains gives an algebraic K-theory invariant of the bundle
called the \emph{higher Franz-Reidemeister (FR) torsion}. Higher FR-torsion distinguishes different smooth structures on the same topological manifold bundle. Therefore, this construction is a strictly differentiable phenomenon.

The purpose of this paper is to explain some of the basic properties
of these constructions and unify them using a simplified version
of Ed Brown's construction. A longer exposition can be found in
\cite{IgComplexTorsion} which, in turn, gives a summary of the
contents of \cite{IgBookOne}.

We summarize the contents of this paper. In Section \ref{sec1: A-infinity functors} we review the definition of an $A_\infty$ functor. $A_\infty$ structures were first constructed by Stasheff \cite{Stasheff63:Polyhedron} and $A_\infty$ categories first appeared in \cite{Fukaya93}. But here we take $A_\infty$ functors only from ordinary categories to the category of $\ZZ$-graded projective modules over a ring $R$. The definition in this restricted case is given by a formula (Equation \ref{eq:def of Ainfty functor}) due to Sugawara\cite{Sugawara60}. We use an old construction of Eilenberg and MacLane\cite{EM53} to make homology into an $A_\infty$ functor when it is projective (Equation \ref {equation of EM}).

In Section \ref{sec2: twisting cochains} we define twisting cochains. We begin with the classical definition of Brown and we also review Brown's
construction of the \emph{twisted tensor product}
\[
    C_\ast(B)\otimes_\f C_\ast(F)
\]
whose homology is equal to the homology of the total space of a
fiber bundle $F\to E\to B$. Brown defined this to be the usual tensor product with boundary map $\d_\f$ twisted by $\f$ (Equation \ref{Brown's twisted boundary}). We give a variation of Brown's definitions which arises from certain $A_\infty$ functors. When there is an underlying functor on the category $\cX$ whose induced maps are all isomorphisms we get a coefficient sheaf $F(X,Y)$ over $\cX$ and the higher homotopies in the $A_\infty$ functor are cochains on $\cX$ with values in $F$. Our twisting cochain is denoted $\psi$ to distinguish it from the classical one of Brown. We use the twisted tensor product to construct a total complex (\ref {def: total complex}) which is an actual functor approximating the $A_\infty$ functor.

Section \ref{sec3:Volodin} describes Volodin K-theory\cite{V} and its
relationship to twisting cochains. We construct the Volodin category $\cV^b(R)$ and a generalization which we call the \emph{Whitehead category} (\ref{def: Whitehead category}). This is a category of acyclic based free chain complexes over a given ring. This category carries a universal twisting cochain of a certain kind and this twisting cochain defines a cohomology class on the classifying space of the Whitehead category:
\[
    \t_{2k}\in H^{4k}(\Wh\bu^h(\QQ,1);\RR).
\]

Higher Franz-Reidemeister torsion is in Section \ref {sec4: higher FR torsion}. We show that, under the right conditions, a smooth fiber bundle with compact
manifold fiber $M\to E\to B$ gives two canonical twisting
cochains, the topologically defined twisting cochain of Ed Brown
and a smoothly defined twisting cochain obtained by fiberwise Morse theory. The
fiberwise mapping cone of the comparison map is fiberwise
contractible. It gives a mapping of the base $B$ into the Whitehead category
provided that we have a basis for the topological twisting
cochain. (The Morse theoretic twisting cochain has a basis coming
from the critical points.) Such a basis can be chosen in the
special case when $\pi_1B$ acts trivially on the rational homology
of $M$. This based free twisting cochain is classified by a map
\[
	B\to |\Wh\bu^h(\QQ,1)|
\]
and we can pull back the universal FR torsion class $\t_{2k}$ to $B$ to obtain the higher FR torsion invariant for the bundle. This invariant has been computed in many cases but here we give only one example: the case when the fiber is a closed oriented even dimensional manifold (Theorem \ref{thm:torsion for M2n closed oriented}).

The rest of the paper contains an elementary discussion of flat
$\ZZ$-grade superconnections. The aim is to show that they are equivalent to twisting cochains. Section \ref{sec5: flat superconnections} derives a definition of an
\emph{infinitesimal twisting cochain}. This is basically a
twisting cochain on very small simplices expressed in terms of
differential forms. The prefix ``super'' refers to a $\ZZ/2\ZZ$ grading. However, a superconnection on a $\ZZ$-graded vector bundle will automatically obtain a $\ZZ$-grading. The Bismut-Lott definition of such a superconnection \cite{Bismut-Lott95} also requires that this $\ZZ$-graded superconnection have total degree 1.

In Section \ref {sec6:supercommutator rules}, we view the endomorphism valued differential forms as operators on the vector bundle in the standard way to obtain flat superconnections. In the last section \ref {sec8: another method} we show that a flat superconnection is the differential in a cochain complex
which is dual to Brown's twisted tensor product. Going backwards, the second to last section \ref{sec7: Chen's iterated integrals} explains how superconnections can be integrated over $1$ and $2$ simplices using Chen's iterated integrals to give the beginning of a simplicial twisting cochain. Complete details for integration of superconnections over arbitrary simplices can be found  in \cite{IgIterated}.

I would like to thank Bernard Keller whose lectures and comments lead me the work of my colleague Ed Brown. I thank Ed Brown for many useful conversations about the topics in this paper. I also thank Jim Stasheff for his support and encouragement and for numerous helpful comments on several versions of this paper and I also want to thank Jonathan Block for explaining his sign conventions to me. This was very helpful for getting the correct signs in \cite{IgIterated} and I have changed the signs in this paper accordingly. Finally, I would like to thank the anonymous referee for numerous suggestions for improving this paper.

%
%

\section{$A_\infty$ functors}\label{sec1: A-infinity functors}

In this paper we consider the differential graded category $\cC(R)$ of chain complexes of projective $R$-modules over an associative ring $R$ and $\cG r(R)$, the underlying graded category of $\ZZ$-graded projective $R$-modules. All $R$-modules will be right $R$-modules. Both categories $\cC(R),\cG r(R)$ have graded hom sets given by
\[
    \HOM(C_\ast,D_\ast)=\bigoplus
    \HOM_n(C_\ast,D_\ast)=\bigoplus_n \prod_k \Hom_R(C_k,D_{n+k}).
\]
But $\cC(R)$ has a differential $m_1:\HOM_n(C_\ast,D_\ast)\to \HOM_{n-1}(C_\ast,D_\ast)$ given by \[m_1(f)=df-(-1)^nfd.\]
We will be considering functors from
an ordinary category $\cX$ into the differential graded category $\cC(R)$, but we view these as functors $\cX\to\cG r(R)$ with additional structure. Usually we assume that the functor takes values in the full subcategory of $\cG r(R)$ consisting either of nonnegatively graded $R$-modules with additional structure given by a degree $-1$ differential or nonpositively graded modules with degree $1$ differential.

To fix a problem which arises in the notation we will use the nerve $\cN\bu\cX\op$ of the opposite category. Thus a \emph{$p$-simplex} in $\cX$ (an element of $\cN_p\cX\op$) will be a
sequence of morphisms of the form:
\begin{equation}\label{eq:p simplex in Xop}
    X_0\xlarrow{f_1}X_1\xlarrow{f_2}\cdots \xlarrow{f_p}X_p
\end{equation}
To clarify the notation, the composition $X_p\to X_0$ is a morphism in $\cX$ which is also a morphism $X_0\to X_p$ in $\cX\op$.
Note that a $0$-simplex consists of one object $X_0$ with no maps. The main purpose of this is to make the domain $X_j$ of the front
$j$-face of a $j+k$ simplex equal to the range of the back
$k$-face.

The following definition of $A_\infty$ functors, in particular Equation (\ref{eq:def of Ainfty functor}), in the restricted case when the domain is an ordinary category is due to Sugawara \cite{Sugawara60}.

\begin{defn}An \emph{$A_\infty$ functor} 
\[
	\Phi=(\Phi,\Phi_0,\Phi_1,\Phi_2,\cdots):\cX\to\cG r(R)
\]
on an ordinary category $\cX$ is an operation which assigns to each object $X\in \cX$ a $\ZZ$-graded projective $R$-module $\Phi X$ and to each sequence of composable morphisms (\ref{eq:p simplex in Xop}) a morphism
\[
    \Phi_p(f_1,f_2,\cdots,f_p):\Phi X_p\to \Phi X_0
\]
of degree $p-1$ and satisfies the following cocycle condition for $p\geq0$.
\begin{multline}\label{eq:def of Ainfty functor}
    \sum_{i=0}^p(-1)^i
    \Phi_i(f_1,\cdots,f_i)\Phi_{p-i}(f_{i+1},\cdots,f_p)\\
    =\sum_{i=1}^{p-1}(-1)^i\Phi_{p-1}(f_1,\cdots,f_if_{i+1},\cdots,f_p)
\end{multline}
For $p\ge1$ this can be written as follows where $m_1(f)=\Phi_0 f-(-1)^{\deg f}f\Phi_0$ and $m_2(f,g)=f\circ g$.
\[
    m_1(\Phi_p)+\sum_{i=1}^{p-1}(-1)^i
    m_2(\Phi_i,\Phi_{p-i})=\sum_{i=1}^{p-1}(-1)^i\Phi_{p-1}(1_{i-1},m_2,1_{p-i-1}).
\]
\end{defn}

For $p=0,1,2,3$ this equation has the following interpretation.
\begin{enumerate}
\item[$(p=0)$] $\Phi_0\Phi_0=0$, i.e., $(\Phi X,\Phi_0)$ is a chain complex.
\item[$(p=1)$]
$\Phi_0\Phi_1(f)=\Phi_1(f)\Phi_0$, i.e., 
\[
	\Phi_1(f):(\Phi X_1,\Phi_0)\to(\Phi X_0,\Phi_0)
\]
is a chain map.
\item[$(p=2)$]
    $\Phi_0\Phi_2(f_1,f_2)+\Phi_2(f_1,f_2)\Phi_0=\Phi_1(f_1)\Phi_1(f_2)-\Phi_1(f_1f_2)$,
    i.e., \[\Phi_2(f_1,f_2):\Phi_1(f_1f_2)\simeq\Phi_1(f_1)\Phi_1(f_2)\] is a chain homotopy.
\item[$(p=3)$] $\Phi_0\Phi_3(f_1,f_2,f_3)-\Phi_3(f_1,f_2,f_3)\Phi_0=$
    \[
    \Phi_2(f_1,f_2f_3)-\Phi_2(f_1,f_2)\Phi_1(f_3)+\Phi_1(f_1)\Phi_2(f_2,f_3)
    -\Phi_2(f_1f_2,f_3).
    \]
    In other words, $\Phi_3$ is a null homotopy of the coboundary
    of $\Phi_2$.
\end{enumerate}

\begin{prop}
Suppose that $\Phi_2=0$ and $\Phi_1$ takes isomorphisms to isomorphisms. Then $(\Phi,\Phi_0,\Phi_1)$ is a functor from the category $\cX$ to the category of projective $R$-complexes and degree $0$ chain maps and $(\Phi,\Phi_1)$ is a functor from $\cX$ to the category of $\ZZ$-graded projective $R$-modules and degree $0$ maps.
\end{prop}

\begin{proof}
If $\Phi_2=0$ then $\Phi_1(f)$ is a chain map $(\Phi X_1,\Phi_0)\to (\Phi X_0,\Phi_0)$ with the property that $\Phi_1(fg)=\Phi_1(f)\Phi_1(g)$. This implies $\Phi_1(id_X)$ is a projection operator. If $\Phi_1$ takes isomorphisms to isomorphisms this must be the identity map on $\Phi X$.
\end{proof}

A \emph{natural transformation} of $A_\infty$ functors is an
$A_\infty$ functor on the product category $\cX\times I$ where $I$
is the category with two objects $0,1$ and one nonidentity
morphism $0\to 1$. This is a family of chain maps $(\Phi X,\Phi_0)\to(\Phi'
X,\Phi_0)$ which are natural only up to a system of higher homotopies. If
these chain maps are homotopy equivalences we say that
$\Phi,\Phi'$ are \emph{$A_\infty$ homotopy equivalent} or
\emph{fiber homotopy equivalent}. (We view an $A_\infty$ functor
as a family of chain complexes over the nerve of $\cX\op$.)

One way to construct an $A_\infty$ functor on $\cX$ is to start
with an actual functor $C$ from $\cX$ to the category of projective
$R$-complexes, then replace each $C(X)$
with a homotopy equivalent projective $R$-complex $\Phi X$ with differential $\Phi_0$.
Following Eilenberg and MacLane \cite{EM53}, the higher
homotopies are given by
\begin{equation}\label{equation of EM}
    \Phi_p(f_1,\cdots,f_p)=q_0C(f_1)\eta_1C(f_2)\eta_2
    \cdots\eta_{p-1}C(f_p)j_p
\end{equation}
where $j_i:\Phi X_i\to C(X_i)$ and $q_i:C(X_i)\to\Phi X_i$ are
homotopy inverse chain maps and $\eta_i:C(X_i)\to C(X_i)$ is a
chain homotopy $id\simeq j_i\circ q_i$.

In the special case when the homology of $C(X)$ is projective
(e.g., if $R$ is a field) and $C(X)$ is either nonnegatively or nonpositively graded, the \emph{homology complex}
\[\Phi X=H_\ast(C(X))\] (with zero boundary map) gives an example of
a homotopy equivalent chain complex. Using the construction above,
we obtain the \emph{$A_\infty$ homology functor}.

Suppose that $p:E\to B$ is a fiber bundle where the base $B$ is
(the geometric realization of) a simplicial complex. Then for each
simplex $\s$ in $B$ the inverse image $E|\s=p^{-1}(\s)$ is
homeomorphic to $F\times\s$ and thus homotopy equivalent to the
fiber $F$. Taking either a cellular chain complex or the total
singular complex, we get a functor from the category $\simp B$ of
simplices in $B$, with inclusions as morphisms, to the category of
augmented chain complexes. When the homology of $F$ is projective,
the construction above gives an $A_\infty$ homology functor on
$\simp B$.

We will see that this $A_\infty$ functor gives a twisting cochain
on $B$ which, by Ed Brown's twisted tensor product construction,
gives a chain complex for the homology of $E$. But first we want to point out that the dual of an $A_\infty$ functor is also an $A_\infty$ functor.

Suppose that $R=K$ is a field. Then we have the degree-wise duality functor on $\cG r(K)$ sending $V=\bigoplus V_n$ to $V^\ast$ where \[
V_n^\ast=\Hom(V_{-n},K).\] Then morphisms of degree $q$ are sent to morphisms of degree $q$ and the order of composition is reversed. So we get a functor $\cG r(K)\to \cG r(K)\op$. Given an $A_\infty$ functor $\Phi:\cX\to \cG r(K)$ we can compose with this duality functor:
\[
	\cX\to \cG r(K)\to \cG r(K)\op
\]
This is the same as a functor $\Phi^\ast:\cX\op\to\cG r(K)$.

\begin{prop}\label{dual of Ainfty is Ainfty}
The composition of an $A_\infty$ functor $\Phi:\cX\to \cG r(K)$ with the degree-wise duality functor on $\cG r(K)$ gives an $A_\infty$ functor $\Phi^\ast$ on $\cX\op$.
\end{prop}

\begin{proof}
This is very straightforward. The interesting point is that there is no change in signs. Apply duality to Equation (\ref{eq:def of Ainfty functor}) and reverse the order of the morphisms $f_i^\ast$ (since they are begin composed in the opposite order in $\cX\op$). We get:
\[
    \sum_{i=0}^p(-1)^i
    \Phi_{p-i}^\ast(f_p^\ast,\cdots, f^\ast_{i+1})\Phi^\ast_i(f_i^\ast,\cdots,f_1^\ast)
    =\sum_{i=1}^{p-1}(-1)^i\Phi_{p-1}^\ast(f_p^\ast,\cdots,f^\ast_{i+1}f_i^\ast,\cdots,f_1^\ast)
\]
Now multiply both sides by $(-1)^p$, replace $i$ by $p-i$ and make the notation change: $g_i=f_{p-i+1}^\ast$ to see that $\Phi^\ast$ satisfies the definition of an $A_\infty$ functor.
\end{proof}

%
%

\section{Twisting cochains}\label{sec2: twisting cochains}

First we review Ed Brown's original construction \cite{Brown59:TwistedTensor}. Suppose that $\LL$ is a commutative ring, $K$ is a nonnegatively graded $\LL$-coalgebra and $A$ is a graded $\LL$-algebra. Then Brown defined a twisting cochain to be a sum $\f=\sum_{p\ge0}\f_p$ of $\LL$-linear maps $\f_p:K_p\to A_{p-1}$,  so that
\begin{enumerate}
\item $\f_0=0$
\item $\d\f_p=\f_{p-1}\d- \sum_{1\le i\le p-1}\f_i\cup'\f_{p-i}$ where
\[
	\f_i\cup'\f_{p-i}=\mu(\f_i\otimes\f_{p-i})\Delta
\]
where $\mu$ is the multiplication in $A$ and $\Delta$ is the comultiplication of $K$. When this expression is evaluated, a sign of $(-1)^i$ is produced by the Koszul sign rule. (See Definition \ref{def of twisting cochain} below.)
\end{enumerate}
Brown also assumed that $A$ has an augmentation, which we do not require.

Given a differential graded $A$-module $M$, Brown defined the \emph{twisted tensor product} 
\[
	K\otimes_\f M
\]
to be the standard graded tensor product over $\LL$ with differential $\d_\f$ given by
\begin{equation}\label{Brown's twisted boundary}
	\d_\f(x\otimes y)=\d x\otimes y+(-1)^{\deg x} x\otimes\d y-\sum_{(x)} (-1)^{\deg x_{(1)}}x_{(1)}\otimes \f(x_{(2)})y
\end{equation}
where we use Sweedler notation $\Delta x=\sum_{(x)}x_{(1)}\otimes x_{(2)}$ \cite{Sweedler}.

Brown then showed that the total singular complex $S_\ast(E)$ (with coefficients in $\LL$) of the total space of a fiber bundle $F\to E\to B$ over a path connected space is homotopy equivalent to the twisted tensor product 
\[
	S'_\ast(B)\otimes_\f C_\ast(F)
\]
where $K=S'_\ast(B)$ is the subcomplex of the total singular complex of $B$ consisting of singular simplices in $B$ with all its vertices at the base point of $B$, $M=C_\ast(F)$ is any free $\LL$-complex homotopy equivalent to the total singular complex of $F$ and $A$ is the differential graded algebra $A=\HOM(M,M)$. (Brown took $A$ to be the total singular complex of the loop space of $B$ which acts on $S_\ast(F)$.)

The cochains $\f_p,p\ge1$ are given as follows. Each 1-simplex $\s$ of $B$ is a loop in $B$ which induces a holonomy $F\to F$. Then $\f_1:=\s_\ast-1$ where $\s_\ast:C_\ast(F)\to C_\ast(F)$ is the induced chain map. Any 2-simplex $\s:\Delta^2\to B$ gives a homotopy between one loop and the composition of the other two loops: $\d_0\s\d_2\s\simeq\d_1\s$. We take $\f_2(\s)$ to be the corresponding chain homotopy $\d_0\s_\ast\circ\d_2\s_\ast\simeq \d_1\s_\ast$. The construction of $\f_p$ for $p\ge3$ is similar and is a special case of the construction given below.

We now consider twisting cochains $\psi=\sum_{p\ge0}\psi_p$ arising from $A_\infty$ functors. In the case of a fiber bundle $p:E\to B$, the idea is that we take simplices in the base $B$ with distinct vertices. For each vertex $v$, $\psi_0$ is the boundary map of the complex $C_\ast(F_v)$. Each 1-simplex gives, not an endomorphism of a single complex $M=C_\ast(F)$, but morphisms between complexes associated to vertices of the simplex. Thus $\HOM(M,M)$ is replaced by a category of graded $R$-modules. We also work over an associative ring $R$. This simply means that the tensor product $S_\ast(B)\otimes_\f M$ is replaced by a direct sum of (shifted) copies of (various) $M$, one for each free generator of $S_\ast(B)$. (See Definition \ref{def: total complex} below). 

Suppose that $(\Phi,\Phi_1):\cX\to \cG r(R)$ is a functor so that $\Phi_1(f)$ is a degree 0 isomorphism for all $f:X\to Y$. In that case, the graded
bifunctor
\begin{equation}\label{eq:F(X,Y)=HOM(Phi X,Phi Y)}
    F(X,Y)=\HOM(\Phi X,\Phi Y)
\end{equation}
gives a locally trivial coefficient system on the category $\cX$.

Since we are using the nerve of the opposite category, a
\emph{$p$-cochain} on $\cX$ with coefficients in a bifunctor $F$ is
a mapping $\psi$ which assigns to each $p$-simplex
\[
    X_\ast=(X_0\xlarrow{f_1} X_1\xlarrow{f_2}\cdots\xlarrow{f_p} X_p)
\]
in $\cX$ an element $\psi(X_\ast)\in F(X_p,X_0)$. The
\emph{coboundary} of $\psi$ is the $p+1$ cochain given by
\begin{multline}
    \delta\psi(X_0,\cdots,X_{p+1})=(f_1)_\ast\psi(X_1,\cdots,X_{p+1})\\
    +\sum_{i=1}^{p}(-1)^i\psi(X_0,\cdots,\widehat{X_i},\cdots,X_{p+1})
    +(-1)^p(f_{p+1})^\ast\psi(X_0,\cdots,X_p).
\end{multline}

\begin{defn}\label{def of twisting cochain} Given a functor $(\Phi,\Phi_1):\cX\to \cG r(R)$ as above, a \emph{twisting cochain} $\psi$ on $\cX$ with coefficients in $(\Phi,\Phi_1)$ is a sum of cochains $\psi=\sum_{p\geq0}\psi_p$ where $\psi_p$ is a $p$-cochain on $\cX$ with coefficients in the degree $p-1$ part $F_{p-1}$ of the $\HOM(\Phi,\Phi)$ bifunctor $F$ of (\ref{eq:F(X,Y)=HOM(Phi X,Phi Y)}) so that the following condition is satisfied.
\[
    \delta\psi=\psi\cup'\psi.
\]
Here $\cup'$ is the cup product using the Koszul sign rule:
\[
    \psi_p\cup'\psi_q(X_0,\cdots,X_{p+q})=(-1)^p\psi_p(X_0,\cdots,X_p)\psi_q(X_p,\cdots,X_{p+q})
\]
since $\psi_q$ has total odd degree.
\end{defn}

To obtain a classical twisting cochain $\f$, we restrict to the case where $\cX$ has a single object $X$, $\Phi_1(f)$ is the identity map on $\Phi X$ for all morphisms $f$ and $R$ is a commutative ring. We take $K$ to be the free $R$-complex of the nerve of $\cX\op$. This is the differential grade $R$-coalgebra $K=C_\ast(\cN\bu \cX\op)$ which in degree $p$ is freely generated by the set of $p$ simplices $X_0\ot \cdots\ot X_p$ in $\cX$. We let $M$ be the projective $R$-complex $M=(\Phi X,\psi_0)$. Then $A=\HOM(M,M)$ is a differential graded $R$-complex. The functor $F$ is trivial and $\delta f=f\d$. So, the equation $\delta\psi=\psi\cup'\psi$ becomes:
\[
	\psi_{p-1}\d=\psi_0\cup'\psi_p+\psi_p\cup'\psi_0+\sum_{i=1}^{p-1}\psi_i\cup'\psi_{p-i}
\]
Since $\psi_0=\d^M$, we have
\[
	\psi_0\cup'\psi_p+\psi_p\cup'\psi_0=\d^M\psi_p+(-1)^p\psi_p\d^M=\d^A\psi_p
\]
Therefore, $\f_0=0$ and $\f_p=\psi_p$ for all $p\ge1$ gives a twisting cochain in the sense of Brown. The referee has pointed out that, if $\Phi_1(f)$ is not always the identity map on $\Phi X$, we can still get a classical twisting cochain by letting $\f_0=0$, $\f_p=\psi_p$ for $p\ge2$ and
\[
	\f_1=\psi_1+\Phi_1-id_M.
\]

By comparison of definitions we have the following.

\begin{prop}
$\psi$ is a twisting cochain on $\cX$ with coefficients in $(\Phi,\Phi_1)$ if and only if
\[(\Phi,\psi_0,\Phi_1+\psi_1,\psi_2,\psi_3,\cdots)\] is an
$A_\infty$ functor.
\end{prop}

Consider again the $A_\infty$ homology functor
\[
    \s\mapsto H_\ast(E_\s)
\]
of a fiber bundle $E\to B$ over a triangulated space $B$. Suppose that the
fiber $F$ has projective homology. In this case $\psi_0$ and
$\psi_1$ are both zero and the higher homotopies $\psi_p$ for
$p\geq2$ are unique up to simplicial homotopy (over $B\times I$). In this case Ed Brown showed that in his {twisted tensor product} $S'_\ast(B)\otimes_\psi S_\ast(F)$, $S_\ast(F)$ can be replaced with $H_\ast(F)$ and $\psi$ is given by the $A_\infty$ homology functor constructed above.

An easy spectral sequence comparison argument shows that we may replace $S'_\ast(B)$ with any homotopy equivalent differential graded coalgebra. We take the cellular complex $C_\ast(B)$ given by a triangulation of $B$. Then the twisted tensor product $C_\ast(B)\otimes_\psi H_\ast(F)$ is the total complex of the usual bicomplex $C_p(B;H_q(F))$ with boundary map modified by the
twisting cochain $\psi$ as follows:
\begin{equation}\label{eq:twisted boundary}
    \d_\psi(x\otimes y)=\d x\otimes y-\sum_{p+q=\deg
    x}(-1)^pf_p(x)\otimes \psi_q(b_q(x))(y).
\end{equation}
Here $x=(x_0\supseteq x_1\supseteq\cdots\supseteq x_n)$ is a
simplex in the first barycentric subdivision of $B$,
$f_p(x)=(x_0\supseteq\cdots\supseteq x_p)$ is the front $p$-face
of $x$ and $b_q(x)=(x_p\supseteq\cdots\supseteq x_n)$ is the back
$q$-face.

\begin{thm}[Brown\cite{Brown59:TwistedTensor}]\label{thm:Brown}
Assuming that $F$ has projective homology, the twisted tensor
product gives the homology of the total space:
\[H_\ast(C_\ast(B)\otimes_\psi H_\ast(F))\cong H_\ast(E).\]
\end{thm}

\begin{rem}\label{rem:psi q gives d in Serre SS}
The mapping $f_p\otimes\psi_qb_q$ in (\ref{eq:twisted boundary})
has bidegree $(-q,q-1)$. It gives the corresponding boundary map
in the Serre spectral sequence for $E$ which is given by filtering
the twisted tensor product by reverse filtration of $H_\ast(F)$
(by subcomplexes $H_{\ast\geq n}(F)$) \cite{IgBookOne}.
\end{rem}

\begin{cor}
When $E\to B$ is an oriented $n-1$ sphere bundle, the degree $n$
part $\psi_n$ of the twisting cochain $\psi$ is a cocycle
representing the Euler class of $E$:
\[
    [\psi_n]=e^E\in H^n(B;R).
\]
\end{cor}

\begin{proof}
Since $\HOM(H_\ast(S^{n-1}),H_\ast(S^{n-1}))$ has elements only in
degrees $0,n-1$, $\psi_k=0$ for $k\noteq n$. By definition of a
twisting cochain we have
\[
    \delta \psi_n=(\psi\cup'\psi)_n=0.
\]
Therefore, $\psi_n$ is an $n-1$ cocycle on $B$. Since it gives the
differential in the Serre spectral sequence, it must represent the
Euler class.
\end{proof}

In the present setting, Ed Brown's twisted tensor product is
equivalent to the following construction.

\begin{defn}\label{def: total complex}
The \emph{total complex} $E(\psi;\Phi)$ of the twisting cochain
$\psi$ with coefficients in the functor $(\Phi,\Phi_1)$ is given
by
\[
    E(\psi;\Phi)=\bigoplus_{k\geq0}\bigoplus_{(X_0\ot\cdots\ot X_k)}
    (X_\ast)\otimes \Phi X_k
\]
with boundary map $\d_\psi$ given by (\ref{eq:twisted boundary}).
\end{defn}

\begin{rem}\label{rem:the total complex of a simplex}
Note that every simplex $X_\ast=(X_0\ot\cdots\ot X_k)$ gives a
subcomplex of the total complex $E(\psi;\Phi)$
by:
\[
    E(X_\ast)=\bigoplus_{j\geq0}\bigoplus_{a:[j]\to[k]}
    a^\ast(X_\ast)\otimes\Phi X_{a(j)}.
\]
This is also the total complex of the $A_\infty$ functor on $[k]$
(considered as a category with objects $0,1,\cdots,k$ and
morphisms $k\to k-1\to\cdots \to1\to0$) given by pulling back
$\psi$ along the functor $X_\ast:[k]\to\cX$.

Note that this gives a functor from the category of simplices in
$\cX$ to the category of subcomplexes of the total complex with
morphisms being inclusion maps. Thus, just as the $A_\infty$
homology functor constructs an $A_\infty$ functor out of an actual
functor, the total complex construction gives an actual functor on
$\simp\cN\bu\cX\op$ from an $A_\infty$ functor on $\cX$.
\end{rem}

Using the total complex, a twisting cochain on $\cX$ can be viewed
as a family of chain complexes parametrized by the nerve of
$\cX\op$. With some extra structure, this gives a map from the
geometric realization of $\cX$ to the Volodin K-theory space of
$R$.

%
%

\section{Volodin K-theory}\label{sec3:Volodin}

Algebraic K-theory is related to twisting cochains in the
following way. When two \emph{based, upper triangular twisting
cochains} are homotopy equivalent, there is an algebraic K-theory
obstruction to deforming one into the other. Formally, we take the
pointwise mapping cone. This gives a based free acyclic upper
triangular twisting cochain on the category $\cX$. This is
equivalent to a mapping from the geometric realization $|\cX|$ of
$\cX$ to a fancy version of the Volodin K-theory space of the ring
$R$. To avoid confusion, we assume that $R$ has the property that the rank of a free $R$-module is well defined, i.e., that $R^n\cong R^m$ implies $n=m$.

When the basis is only well-defined up to permutation and
multiplication by elements of a subgroup $G$ of the group of units
of $R$, an acyclic twisting cochain on $\cX$ defines a mapping
from $|\cX|$ into the fiber $\Wh\bu^h(R,G)$ of the mapping
\[
    \Omega^\infty\Sig^\infty(BG_+)\to BGL(R)^+\times\ZZ.
\]
The well-known basic case is the Whitehead group
\[
    Wh_1(G)=\pi_0 \Wh\bu^h(\ZZ[G],G)
\]
which is the obstruction to $G$-collapse of a contractible f.g.
based free $R$-complex. In this section we discuss the different
versions of the Volodin construction, show how they are related to
twisting cochains and identify the homotopy type of two of them.

The basic definition is sometimes called the ``one index'' case.
It is a space of invertible matrices locally varying by upper
triangular column operations. When this definition is expressed as
a twisting cochain, the construction seems artificial, with only
one term $\psi_1$ in the twisting cochain $\psi$. However, when
the missing higher terms $\psi_p$ are inserted we recover the
general Volodin space.

\begin{defn} For every $n\geq 2$ the \emph{Volodin category} $\cV^n(R)$ is the category whose objects are pairs $(A,\s)$
consisting of an invertible $n\times n$ matrix $A\in GL(n,R)$ and
a partial ordering $\s$ of $\{1,2,\cdots,n\}$. A morphism
$(A,\s)\to(B,\t)$ is an $n\times n$ matrix $T$ with coefficients in $R$ so that
\begin{enumerate}
    \item $\s\subseteq\t$, i.e., $\t$ is a refinement of $\s$.
    \item $AT=B$. (So the morphism $T=A^{-1}B$ is unique if it exists.)
    \item $T=(t_{ij})$ is \emph{$\t$-upper triangular} in the sense that\begin{enumerate}
        \item $t_{ii}=1$ for $i=1,\cdots,n$
        \item $t_{ij}=0$ unless $i\leq j$ in the partial ordering
    $\t$ ($i\le j\Leftrightarrow (i,j)\in\t$).
    \end{enumerate}
\end{enumerate}
Note that composition is reverse matrix multiplication:
\[
    S\circ T=TS.
\]
\end{defn}

There is a \emph{simplicial Volodin space} $V\bu^n(R)$ without
explicit partial orderings. A $p$-simplex $g\in V_p^n(R)$ consists
of a $p+1$ tuple of invertible $n\times n$ matrices
\[
    g=(g_0,g_1,\cdots,g_p)
\]
so that for some partial ordering $\s$ of $\{1,\cdots,n\}$ the
matrices $g_i^{-1}g_j$ are all $\s$-upper triangular. There is a
simplicial map
\begin{equation}\label{eq:N Vn(R) to Vn(R)}
    \cN\bu\cV^n(R)\to V\bu^n(R)
\end{equation}
from the nerve of the Volodin category $\cV^n(R)$ to the simplicial
set $V\bu^n(R)$ given by forgetting the partial orderings.
However, the collection of admissible partial orderings on any
$g\in V_p^n(R)$ has a unique minimal element and therefore forms a
contractible category. Consequently, (\ref{eq:N Vn(R) to Vn(R)})
induces a homotopy equivalence
\[
    |\cV^n(R)|\simeq |V\bu^n(R)|
\]

If we stabilize matrices in the usual way by adding a $1$ in the
lower right corner we get the stable Volodin category
\[
    \cV(R)=\lim_\to\cV^n(R)
\]
and the stable Volodin space $\displaystyle{V\bu^\infty(R)=\lim_\to V\bu^n(R)}$ which are related to Quillen K-theory by the following well-known
theorem due to Vasserstein and Wagoner but best explained by
Suslin \cite{Suslin:equivalence-of-K-theories}.

\begin{thm}\label{thm:equivalence of K-theories}
$|\cV(R)|\simeq|V\bu^\infty(R)|\simeq \Omega BGL(R)^+$ where
$\Omega BGL(R)^+$ is the loop space of the plus construction on
the classifying space of $GL(R)=GL(\infty,R)$.
\end{thm}

The Volodin category $\cV^n(R)$ has a canonical twisting cochain.
It comes from the realization that an invertible matrix is the
same as a based contractible chain complex with two terms. 

\begin{defn}\label{def: canonical twisting cochain on V(R)} The
\emph{canonical twisting cochain} on $\cV^n(R)$ is given as follows.
\begin{enumerate}
    \item Let $\Phi(A,\s)=C_\ast$ be the based free graded $R$-module with $C_0=C_1=R^n$ for every object $(A,\s)$ of
    $\cV^n(R)$.
    \item $\Phi_1=(id,id)$ is the identity chain map $C_\ast\to
    C_\ast$ for every morphism.
    \item $\psi_0(A,\s)=A:R^n\to R^n$.
    \item $\psi_1(T)=(0,T-I)$. So $I+\psi_1(T)=(I,T)$ gives a
    chain isomorphism:
\begin{displaymath}
\begin{CD}
R^n  @>T>>  R^n\\
@V{B}VV               @VV{A}V \\
R^n   @>I>>    R^n
\end{CD}
\end{displaymath}
\end{enumerate}
\end{defn}

The higher homotopies $\psi_p, p\geq2$, are all zero for
$\cV^n(R)$. However, there is a fancier version of the Volodin
category with higher homotopies. We call it the ``Whitehead
category.'' This is very similar to the original definition of Volodin \cite{V}

\begin{defn}\label{def: Whitehead category}
If $G$ is a subgroup of the group of units of a ring $R$ then the
\emph{Whitehead category} $\Wh\bu(R,G)$ is defined to be the
simplicial category whose simplicial set of objects consists of
pairs $(P,\psi)$ where $\psi$ is an upper-triangular twisting
cochain on the category $[k]$ (as in Remark \ref{rem:the total
complex of a simplex}) with coefficients in the fixed graded based $R$-module:
\[
    R^P:=\bigoplus R^{P_i}
\]
where $P=\coprod P_i$ is a graded poset (a poset with a grading
not necessarily related to the ordering). By
\emph{upper-triangular} we mean that $\psi(\s)(x)$ is a linear
combination of $y<x\in P$ for all simplices $\s$ in $[k]$. As in
the Volodin category, $\Phi_1$ is the identity mapping on $R^P$.

A \emph{morphism} $(P,\psi)\to(Q,\f)$ in $\Wh_k(R,G)$ consists of
a graded order preserving monomorphism going the other way:
\[
    f:Q\to P
\]
so that $S=P-f(Q)$ is a disjoint union of \emph{expansion pairs}
which are, by definition, pairs $x^+>x^-$ otherwise unrelated to
every other element of $P$ so that $\psi_0(x^+)=x^-g$ for some
$g\in G$, together with a function $\g:Q\to G$ so that $\f$
differs from $\psi\circ f$ only by multiplication by $\g$, i.e.,
$f_\ast(\f_p(\s)(x))=\psi_p(\s)(f(x))\g(x)$.
\end{defn}

The following theorem, due to J. Klein and the author, is proved
in \cite{IgBookOne}, Section 5.6.

\begin{thm}[Igusa-Klein]\label{thm:homotopy type of Whitehead
category} There is a homotopy fiber sequence
\[
    |\Wh\bu^h(R,G)|\to \Omega^\infty \Sig^\infty(BG_+)\to \ZZ\times BGL(R)^+
\]
where $\Wh\bu^h(R,G)$ is the simplicial full subcategory of
$\Wh\bu(R,G)$ consisting of $(P,\psi)$ so that each chain complex
$(R^P,\psi_0)$ is contractible (i.e., has the homology of the
empty set) and $BG_+=BG\coprod pt$.
\end{thm}

\begin{rem}
If we take $G$ to be a finite group then $\Omega^\infty
S^\infty(BG_+)$ is rationally trivial above degree $0$ so
$\Wh\bu^h(R,G)$ has the rational homotopy type of the Volodin
space:
\[
    |\Wh\bu^h(R,G)|\simeq_\QQ |\cV(R)|\simeq \Omega BGL(R)^+.
\]
In particular, if $R=\QQ$, we get \cite{Borel74}
\[
    |\Wh\bu^h(\QQ,1)|\simeq_\QQ \Omega BGL(\QQ)^+\simeq_\QQ BO.
\]
\end{rem}

Using the Borel regulator maps
\[
    K_{4k+1}\QQ\to\RR
\]
given by continuous cohomology classes in $H^{2k}(BGL(\CC);\RR)$
we get the \emph{universal real higher Franz-Reidemeister torsion}
invariants
\[
    \t_{2k}\in H^{4k}(\Wh\bu^h(\QQ,1);\RR).
\]
These give characteristic classes for smooth bundles under certain
conditions.

%
%

\section{Higher FR torsion}\label{sec4: higher FR torsion}

We will discuss the circumstances under which we obtain well
defined algebraic K-theory classes for a fiber bundle. If we have
a smooth bundle $p:E\to B$ where $E, B$ and the fiber
$M_b=p^{-1}(b)$ are compact connected smooth manifolds and $R$ is
a commutative ring so that the fiber homology $H_\ast(M_b;R)$ is
projective then we obtain two canonical twisting cochains on $B$.

The first is Brown's twisting cochain $\psi$ with coefficients in
the fiberwise homology bundle \[\Phi(b)=H_\ast(M_b;R).\] Recall
that this requires the fiber homology to be projective.

The second is the fiberwise cellular chain complex $C_\ast(f_b)$
associated to a fiberwise \emph{generalized Morse function} (GMF)
$f:E\to\RR$. These are defined to be smooth functions which, on
each fiber $M_b$, have only Morse and birth-death singularities
(cubic in one variable plus nondegenerate quadratic in the
others). The fiberwise GMF is \emph{not} well-defined up to
homotopy. However, there is a canonical choice called a ``framed
function'' (\cite{IgFF}, \cite{IgBookOne},
\cite{IgComplexTorsion}) which exists stably and is unique up to
framed fiber homotopy. This gives the following.

\begin{thm}\label{thm:cellular chain complex of a fiberwise framed
function} Any compact smooth manifold bundle $E\to B$ gives a
mapping
\[
    C(f): B\to|\Wh\bu(\ZZ[\pi_1E],\pi_1E)|
\]
which is well-defined up to homotopy and fiber homotopy equivalent
to the fiberwise total singular complex of $E$ with coefficients
in $\ZZ[\pi_1E]$.
\end{thm}

\begin{rem}
The \emph{fiberwise total singular complex} of $E$ is the functor
which assigns to each simplex $\s:\Delta^k\to B$, the total
singular complex of $E|\s$. The fiberwise framed function $f$ is
defined on a product space $E\times D^N$.
\end{rem}

In order to compare the two constructions we need a representation
\[
    \r:\pi_1E\to U(R)
\]
of $\pi_1E$ into the group of units of $R$ with respect to which the fiber homology $H_\ast(M_b;R)$ is projective over $R$. By the functorial properties of the Whitehead category we get a mapping
\[
    B\to|\Wh\bu(R,G)|
\]
where $G\subseteq U(R)$ is the image of $\r$. By Theorem
\ref{thm:cellular chain complex of a fiberwise framed function},
this will be fiberwise homotopy equivalent to the $A_\infty$
fiberwise homology functor $\Phi_1+\psi$. The fiberwise mapping
cone will be fiberwise contractible but it will not give a mapping
to $\Wh\bu^h(R,G)$ unless the fiberwise homology has a basis. This
gives the following.

\begin{cor}\label{cor:existence of higher torsion}
If $\pi_1B$ acts trivially on the fiberwise homology
$H_\ast(M_b;R)$ then a fiberwise mapping cone construction gives a
mapping
\[
    C(C(f)):B\to |\Wh\bu^h(R,G)|
\]
which is well-defined up to homotopy.
\end{cor}

\begin{rem}
A more precise statement is that we take the direct sum of the
fiberwise mapping cone with a fixed contractible projective
$R$-complex $P_\ast$ with the property that $H_\ast(M_b;R)\oplus
P_\ast$ is free in every degree.
\end{rem}

When $R=\QQ$, the construction of higher FR torsion extends to the case then $\pi_1B$ acts \emph{unipotently} on $H_\ast(M;\QQ)$ by which we mean that $H_\ast(M;\QQ)$ admits a filtration by $\pi_1B$ submodules so that the action of $\pi_1B$ on the successive subquotients is trivial.

\begin{cor}
Suppose that $E\to B$ is a compact smooth manifold bundle over a
connected space $B$ so that $\pi_1B$ acts unipotently on the
rational homology of the fiber $M$. Then we have a mapping
\[
    B\to |\Wh\bu^h(\QQ,1)|
\]
which is well-defined up to homotopy and we can pull back
universal higher torsion invariants to obtain well-defined
cohomology classes
\[
    \t_{2k}(E)\in H^{4k}(B;\RR)
\]
which are trivial if the bundle is diffeomorphic to a product
bundle.
\end{cor}

It has been known for many years (by
\cite{Farrell-Hsiang:Diff(Dn)} using the stability theorem
\cite{IgStability}) that there are smooth bundles which are
homeomorphic but not diffeomorphic to product bundles and that
these exotic smooth structures are detected by algebraic K-theory.
Therefore, when the higher FR torsion was successfully defined, it
had already been known to be nonzero in these cases.

However, in these exotic examples the fiber $M$ is either odd
dimensional or even dimensional with boundary. We now have the complete calculation of the higher torsion in the
case of closed oriented even dimensional fibers.

\begin{thm}[6.6 in \cite{IgComplexTorsion}]\label{thm:torsion
for M2n closed oriented} Suppose that $M^{2n}$ is a closed
oriented even dimensional manifold and $M\to E\to B$ is a smooth
bundle so that $\pi_1B$ acts unipotently on the rational homology of
$M$. Then the higher FR torsion invariants $\t_{2k}(E)$ are
well-defined and given by
\[
    \t_{2k}(E)=\frac12(-1)^k\z(2k+1)\frac1{(2k)!}T_{2k}(E)\in
    H^{4k}(B;\RR)
\]
where $\z(s)=\sum \frac1{n^s}$ is the Riemann zeta function and
\[
    T_{2k}(E)=tr_B^E\left(\frac{(2k)!}{2}ch_{2k}(T^vE\otimes\CC)\right)\in
    H^{4k}(B;\ZZ)
\]
with $T^vE$ being the vertical tangent bundle of $E$, $ch_{2k}(T^vE\otimes\CC)$ stands for the degree $4k$ term in the Chern character of $T^vE\otimes\CC$ and
\[
    tr_B^E:H^{n}(E;\ZZ)\to H^n(B;\ZZ)
\]
is the transfer (with $n=4k$).
\end{thm}

\begin{rem}
Note that $T_{2k}(E)$ is a tangential fiber homotopy invariant.
This is in keeping with the belief that there are rationally no stable exotic smooth structures on bundles with closed oriented even dimensional fibers. (\emph{Stable} means stable under product with large dimensional
disks $D^N$. The exotic smooth structure on disk bundles and odd
dimensional sphere bundles of \cite{Farrell-Hsiang:Diff(Dn)} and
the explicit examples given by Hatcher (\cite{IgBookOne},\cite{Goette01}) are stable.) For more details about this conjecture see \cite{GoetteIgusa10}. In that paper we construct virtually all stable exotic smooth structures on bundles with closed odd dimensional fibers and explain why the even dimensional case is so different. See also \cite{IgPontjagin} and  \cite{Goette08} for an outline of those results.

In the special case when $n=1$, $M$ is an oriented surface and the
bundle $E$ is classified by a map of $B$ into the classifying
space $BT_g$ of the Torelli group $T_g$ where $g$ is the genus of
$M$. The tangential homotopy invariant $T_{2k}$ is equal to the
Miller-Morita-Mumford class in this case
(\cite{Mumford83:MMM_class}, \cite{Morita84},
\cite{Miller86:MMM}).  It is still unknown whether or not any of
these classes (tautological classes in degree $4k$) is rationally nontrivial on the Torelli group.
\end{rem}

There are several competing versions of higher FR torsion and the version described here is sometimes called \emph{Igusa-Klein (IK) torsion} since the first computation was given in \cite{IK1Borel2}. Dwyer, Weiss and Williams have defined three kinds of higher Reidemeister which are called \emph{smooth, topological and homotopy DWW torsion} \cite{DWW}. Badzioch, Dorabiala, Klein and Williams have recently shown \cite{BDKW} that smooth DWW torsion is equivalent to IK torsion, making use of the axiomatic characterization of higher torsion given in \cite{IAxioms0}.

Bismut and Lott \cite{Bismut-Lott95} have defined higher
analytic torsion invariants which have been computed in many cases
(\cite{Bismut-Lott97}, \cite{BG2}, \cite{Bunke:spheres},
\cite{Ma97}, \cite{Goette01}). In the case of closed oriented
even dimensional fibers, the analytic torsion is always zero and
Goette has now shown that the expression in Theorem
\ref{thm:torsion for M2n closed oriented} gives the difference
between BL torsion and IK torsion in all cases. (See the survey article \cite{Goette08}.)

Higher analytic torsion is defined using flat $\ZZ$-graded
superconnections. It was observed by Goette \cite{Goette01} that
these are infinitesimal twisting cochains or, as he puts it, that
twisting cochains are combinatorial superconnections. We will
explain this comment.

%
%

\section{Flat superconnections}\label{sec5: flat superconnections}

When we review the definition of a flat $\ZZ$-graded
superconnection, we will see that it is the same as an
infinitesimal twisting cochain. More precisely, the
superconnection is the boundary map of the infinitesimal twisted
tensor product. This gives one explanation of the supercommutator
rules.

Instead of defining superconnections and showing their
relationship to twisting cochains we will take the opposite
approach. We ask the question: What is the natural definition of
an ``infinitesimal twisting cochain?'' This question will lead us
to the definition of a flat superconnection and we will see that
the ``superconnection complex'' $(\Omega(B,V),D)$ is dual to a
twisted tensor product.

Suppose that $B$ is a smooth manifold and $C=\bigoplus_{n\geq0}
C_n$ is a nonnegatively graded complex vector bundle over $B$. Suppose we have a graded flat connection $\nabla$ on $C$ making each $C_n$ into a locally constant coefficient sheaf for the twisting cochain
that we want. The example that we keep in mind is when $C$ is the
fiberwise homology of a smooth manifold bundle $F\to E\xrarrow p B$. By this we mean the graded vector bundle over $B$ whose fiber over $b\in B$ is the homology of $p^{-1}(b)$. The dual bundle
\[
    C^\ast:=\bigoplus_{n\geq0}\Hom(C_n,\CC)
\]
is the fiberwise cohomology bundle $H^\ast(F)\to C^\ast\to B$.

Now, imagine that $B$ is subdivided into tiny simplices and we
have a twisting cochain on $B$ with coefficients in $(C,\nabla)$ which satisfies smoothness conditions to be added later.
Then, at each vertex $v$ we have a degree $-1$ endomorphism
$\psi_0(v)$ of $C(v)$. This gives a degree $1$ endomorphism
$A_0=\psi_0^\ast$ of the dual $C^\ast(v)$. Suppose we can extend this
to a smooth family of such maps
\[
    A_0\in\Omega^0(B,\End(C^\ast))=\Omega^0(B,\End(C)\op)
\]
so that $A_0(x)$ has degree $1$ and square zero ($A_0(x)^2=0$) at
all $x\in B$.

Next, we take the edges of $B$. If an edge $e$ goes from $v_0$ to
$v_1$ the twisting cochain gives us a degree $0$ map
\[
    C(v_0)\xlarrow{\psi_1(e)} C(v_1)
\]
so that $\psi_1(e)$ together with the map (parallel transport)
given by the flat connection $\nabla$ is a chain map. This chain map is the parallel transport of a non-flat connection $\nabla_1$ which we now describe.

If we dualize $\psi_1(e)$ and take only the linear term (ignoring $\Delta v^2$ terms) we get a degree $0$ map $A_1(\Delta v):C^\ast(v_0)\to C^\ast(v_1)$ which is linear in $\Delta v$. To obtain the smooth version we need to take local coordinates for $C$ so that parallel transport of $\nabla$ is constant, i.e., so that, on $C^\ast$, $\nabla^\ast=d$. Then $A_1$ becomes a matrix $1$-form on $B$ (assuming the twisting cochain is smooth in a suitable sense)
\[
    A_1\in\Omega^1(B,\End(C^\ast))
\]
so that parallel transport by the new connection $\nabla_1=d-A_1$
on $C^\ast$ keeps $A_0$ invariant. (The change in sign comes from
the fact that parallel transport by $d-A_1$ is given
infinitesimally by $I+A_1(\Delta v)$ where $A_1(\Delta v)$ is evaluation of the matrix 1-form $A_1$ on the vector $\Delta v$.) This means that
\[
    [\nabla_1,A_0]=[d-A_1,A_0]=0
\]
\[
    dA_0=[A_1,A_0]=A_1A_0+A_0A_1.
\]


We interpret this as an approximately commutative diagram:
\begin{displaymath}
\begin{CD}
    C^\ast(v_0)  @>I+A_1(\Delta v)>>  C^\ast(v_1)\\
    @A{A_0}AA               @AA{A_0+\Delta A_0}A \\
    C^\ast(v_0)   @>I+A_1(\Delta v)>>    C^\ast(v_1)
\end{CD}
\end{displaymath}
Higher order terms are needed to make the diagram actually commute. The linear terms give the following approximate equation:
\[
	\Delta A_0\cong A_1(\Delta v)A_0-A_0A_1(\Delta v)
\]
Since $A_0$ is odd, we get two changes of signs:
\[
	\Delta A_0\cong -dA_0(\Delta v)\quad\quad A_1(\Delta v)A_0=-(A_1A_0)(\Delta v)
\]
As $\Delta v\to0$ we get the equation $dA_0=A_1A_0+A_0A_1$ as claimed.

At the next step, we take two small triangles in $B$ forming a
rectangle. The following diagram which commutes up to homotopy by
$\psi_2^\ast=A_2$ indicates what is happening. Here
$A_1=A_1^xdx+A_1^ydy$ where $A_1^x,A_1^y$ are (even) matrix $0$-forms and $A_1^x \Delta x$ indicates multiplication by the scalar quantity $\Delta x$.
\begin{displaymath}
\begin{CD}
    C^\ast(v_1')  @>{I+A_1^x\Delta x+ \Delta A_1^x \Delta x}>>  C^\ast(v_2)\\
    @A{I+A_1^y \Delta y}AA               @AA{I+A_1^y \Delta y+ \Delta A_1^y \Delta y}A \\
    C^\ast(v_0)   @>\quad\  I+A_1^x \Delta x\quad\ >>    C^\ast(v_1)
\end{CD}
\end{displaymath}
This gives the following approximate equation where $\Delta x, \Delta y$ are scalar quantities and $\Delta v_x,\Delta v_y$ are the corresponding vector quantities giving our rectangle in $B$.
\begin{multline*}
    (A_0A_2+A_2A_0)(\Delta v_x,\Delta v_y)\cong A_1^y \Delta y+ \Delta A_1^y \Delta y A_1^x \Delta x
    -A_1^x \Delta x+ \Delta A_1^x \Delta x A_1^y \Delta y\\
    \cong A_1^yA_1^x \Delta x \Delta y+\frac{\d A_1^y}{\d x} \Delta x \Delta y
    -A_1^xA_1^y \Delta x \Delta y
    -\frac{\d A_1^x}{\d y} \Delta x \Delta y
\end{multline*}
Since $A_1^x,A_1^y$ are even, the right hand side can
be written as
\[
    \left(dA_1-A_1^2\right)(\Delta v_x,\Delta v_y)
\]
In other words, we have
\[
    A_2\in\Omega^2(B,\End(C^\ast))
\]
satisfying the equation
\[
    dA_1=A_0A_2+A_1^2+A_2A_0.
\]

In general we will require that 
\[
    dA_{n-1}=\sum_{p+q=n}A_pA_q.
\]
(See \cite{IgIterated} for a full explanation.) This leads to the following
definition.

\begin{defn}
An \emph{infinitesimal twisting cochain} on $B$ with coefficients
in a graded vector bundle $C^\ast$ with graded flat connection $\nabla$ ($\nabla=\sum \nabla_k$ where $(-1)^k\nabla_k$ is a flat connection on $C^k$)
is equal to a sequence of $\End(C^\ast)$-valued forms
\[
    A_p\in\Omega^p(B,\End_{1-p}(C^\ast))=\Omega^0(B,\End_{1-p}(C^\ast))\otimes_{\Omega^0(B)}\Omega^p(B)
\]
of total degree $1$ so that
\begin{equation}\label{eq:infinitesmal twisting cochain}
    \nabla A_{n-1}=\sum_{p+q=n}A_pA_q.
\end{equation}
\end{defn}

Next we pass to the algebra of operators on $\Omega(B,C^\ast)$ where we carefully distinguish between differential forms $A$ and the operators $\widetilde A$ that they define to arrive at the Bismut-Lott definition of a flat $\ZZ$-grade superconnection.

%
%

\section{Forms as operators}\label{sec6:supercommutator rules}

If $A\in\Omega(B,\End(C^\ast))$ is written as $A=\sum \f_i\otimes
\a_i$ with fixed total degree $|A|=|\f_i|+|\a_i|$, let
$\widetilde{A}$ be the linear operator on
\[
    \Omega(B,C^\ast)=\Omega^0(B,C^\ast)\otimes_{\Omega^0(B)} \Omega(B)
\]
given by
\begin{equation}\label{eq:def of widetilde(A)}
    \widetilde{A}(c\otimes \g):=\sum_{i}(-1)^{|c|\cdot|\a_i|} \f_i(c)\otimes \a_i\wedge \g
\end{equation}

\begin{prop}[Prop. 1 in \cite{Quillen:superconnections}]\label{prop:widetilde(A) is as operator}
If $\omega\in\Omega^k(B)$ then
\[
    \widetilde{A}\circ\omega=(-1)^{k|A|}\omega\circ\widetilde{A}.
\]
Conversely, any linear operator on $\Omega(B,C^\ast)$ of fixed
total degree having this property is equal to $\widetilde{A}$ for
a unique $A\in \Omega(B,\End(C^\ast))$.
\end{prop}

\begin{proof}
Since $\widetilde A$ acts only on the first tensor factor we get
\[
	\widetilde A\circ \w= \widetilde{A\w}= (-1)^{k|A|}\widetilde{\w A}=(-1)^{k|A|}\w\circ\widetilde{A}
\]
as required. Conversely, any linear operator which is $\Omega(B)$-linear in this sense must be ``local'' and thus we may restrict to a
coordinate chart $U$ over which $C^\ast$ has a basis of sections.
This makes $\Omega(U,C^\ast|U)$ into a free module over
$\Omega(U)$. Thus any $\Omega(U)$-linear operator is uniquely
given by $\widetilde{A}$ where $A\in\Omega(U,\End(C^\ast|U))$ is
given by the value of the operator on the basis of sections of
$C^\ast|U$. We can patch these together on intersections of coordinate charts since, by uniqueness, the differential forms defined using different coordinate charts will agree.
\end{proof}

Here is another straightforward calculation.
\begin{prop}
$[d, \widetilde A]=d\circ \widetilde A-(-1)^{|A|} \widetilde A\circ d= \widetilde{dA}$.
\end{prop}

If $A'$ is another $\End(C^\ast)$ valued
form on $B$ then $\widetilde{AA'}= \widetilde A\circ \widetilde{A'}$. So
\[
    [d,\widetilde{A}_{n-1}]
    =\widetilde{dA}_{n-1}=\sum_{p+q=n}\widetilde{A}_p\circ\widetilde{A}_q
\]
which, in coordinate free notation, is
\[
    [\nabla,\widetilde{A}]=\widetilde{A}\circ\widetilde{A}
\]
Since $|A|=1$ and $\nabla^2=0$, we get
\[
    (\nabla-\widetilde{A})^2=(\nabla-\widetilde{A})\circ(\nabla-\widetilde{A})=0.
\]
This leads to the following definition due to Bismut and Lott
\cite{Bismut-Lott95}. (A similar definition appeared in \cite{Chen75}.)

\begin{defn}
Let $V=\bigoplus_{n\geq0}V^n$ be a graded complex vector bundle
over a smooth manifold $B$. Then a \emph{superconnection} on $V$
is defined to be a linear operator $D$ on $\Omega(B,V)$ of total
degree $1$ so that
\[
    D\a=d\a+(-1)^{|\a|}\a D
\]
for all $\a\in\Omega(B)$. The superconnection $D$ is called
\emph{flat} if
\[
    D^2=0.
\]
\end{defn}

If $D$ is flat then $(\Omega(B,V),D)$ is a chain complex which we call the \emph{superconnection complex}. We will see later that it is homotopy equivalent to the dual of a twisted tensor product. The superconnection complex is bigraded: 
\[
\Omega(B,V)=\bigoplus\Omega^p(B,V^q)
\]and the superconnection $D$ has terms of degree $(k,1-k)$ for $k\ge0$. This gives a spectral sequence in the usual way with $E_1^{p,q}=\Omega^p(B,H^q(V,A_0))$ and
\[
	E_2^{p,q}=H^p(B;H^q(V,A_0))\Rightarrow H^{p+q} (\Omega(B,V),D)
\]

A flat superconnections on $V$ corresponds to a contravariant $A_\infty$
functor on $B$. To get a twisting cochain we need an ordinary
graded flat connection $\nabla$ on $V$. Then $D-\nabla$ gives a
twisting cochain by reversing the above process.

The first step is to get out of the superalgebra framework by
writing $D$ as a sum
\[
    D=\nabla-\widetilde{A}=\nabla-\widetilde{A}_0-\widetilde{A}_1
    -\widetilde{A}_2
    -\cdots
\]
where $A_p\in \Omega^p(B,\End_{1-p}(V))$ corresponds to
$\widetilde{A}_p$ by (\ref{eq:def of widetilde(A)}) and satisfies
(\ref{eq:infinitesmal twisting cochain}).

Next, we obtain a contravariant twisting cochain on the category of smooth simplices in $B$ with coefficients in the category of cochain complexes by
iterated integration of $A_\ast$. Then we dualize, relying on Proposition \ref{dual of Ainfty is Ainfty} to recover the original twisting cochain.

%
%

\section{Chen's iterated integrals}\label{sec7: Chen's iterated integrals}

This section gives a very short discussion and proof of the first two steps in the process of integrating a flat superconnection to obtain a twisting cochain. Details are fully explained in \cite{IgIterated} although the original idea is contained in Chen's work \cite{Chen73}, \cite{Chen75}, \cite{Chen:Iterated-path-integral}. A more direct, less computational method of constructing the twisting cochain is explained in the next section.

Since $\nabla$ is a flat connection, we can choose local
coordinates so that $V$ is a trivial bundle and $\nabla=d$. Starting with $p=0$ we note that $A_0(x)$ is a degree $1$
endomorphism of $V_x$ with $A_0(x)^2=0$ making $C(x)=(V_x,A_0(x))$ into a cochain complex for all $x\in B$. Putting $n=2$ in (\ref{eq:infinitesmal twisting cochain}) we see that the curvature $(d-A_1)^2$ of the connection $d-A_1$ is null homotopic. Also we will see that
parallel transport of this connection is a cochain map.

It is well-known that the parallel transport $\Phi_1$ associated
to the connection $d-A_1$ on $V$ is given by an iterated
integral of the matrix $1$-form $A_1$. Given any piecewise
smooth path $\g:[0,1]\to B$, parallel transport is the family of
degree zero homomorphisms $\Phi_1(t,s):C(\g(s))\to C(\g(t))$ so
that $\Phi_1(s,s)=I=id_V$ and $d-A_1=0$, i.e.,
\[
    \frac{\d}{\d t}\Phi_1(t,s)=A_1/t\Phi_1(t,s)
\]
\[
    \frac{\d}{\d s}\Phi_1(t,s)=-\Phi_1(t,s)A_1/s
\]
where $A_1/t=A_1(\g(t))(\g')\in \End(C(\g(t)))$. The
solution is given by Chen's iterated integral
\cite{Chen:Iterated-path-integral}:
\begin{multline*}
    \Phi_1(s_0,s_1)=I+\int_{s_0\ge t\ge s_1}dt_1A_1/t
    +\int_{s_0\ge t_1\ge t_2\ge
    s_1}dt_1dt_2(A_1/t_1)(A_1/t_2)\\
    +\int_{s_0\ge t_1\ge t_2\ge t_3\ge s_1}dt_1dt_2dt_3 (A_1/t_1)(A_1/t_2)(A_1/t_3)
    +\cdots
\end{multline*}
which we abbreviate as:
\[
    \Phi_1(s_0,s_1)=I+\int_\g A_1+\int_\g
    (A_1)^2+\int_\g
    (A_1)^3+\cdots.
\]
This can also be written as a limit of products (multiplied right to left)
\[
	\Phi_1=\lim_{\Delta t\to 0}\prod (I+(A_1/t_i)\Delta t)
\]
In the case when $A_1$ is constant, parallel transport
$C(\g(0))\to C(\g(1))$ is given by $e^{A_1}$. The inverse is given by
$\Phi_1(0,1)=e^{-A_1}$.

\begin{prop}\label{lemma 1}
$A_0(\g(t))\Phi_1(t,s)=\Phi_1(t,s)A_0(\g(s))$, i.e.,
$\Phi_1(t,s)$ gives a cochain map \[C(\g(t))\ot C(\g(s)).\]
\end{prop}

\begin{proof} By (\ref{eq:infinitesmal twisting cochain}) we have:
\[
   - \frac{d}{dt}A_0(\g(t))=dA_0(\g(t))(\g')
    =A_0(\g(t))A_1/t-(A_1/t)A_0(\g(t))
\] where both negative signs come from the fact that $A_0$ is odd.
So, $X(t)=A_0(\g(t))\Phi_1(t,s)$ is the unique solution of
the differential equation
\[
    \frac{\d}{\d t}X(t)=(A_1/t)X(t)
\]
with initial condition $X(s)=A_0(\g(s))$. So, $X(t)$ must also be equal to the other solution of this differential equation which is $X(t)=\Phi_1(t,s)A_0(\g(s))$.
\end{proof}

Let
\[
    \Delta^2=\{(x,y)\in\RR^2\st 1\geq x\geq y\geq0\}
\]
and suppose that $\s:\Delta^2\to B$ is a smooth simplex with
vertices $v_0=\s(0,0),v_1=\s(1,0),v_2=\s(1,1)\in B$. Then a chain
homotopy
\[
    \Phi_1(v_0,v_2)\simeq\Phi_1(v_0,v_1)\Phi_1(v_1,v_2)
\]
can be obtained by an iterated integral of the form
\[
    \psi_2(\s)=\int_\s A_2+\int_\s A_2 A_1+\int_\s
    A_1 A_2+\int_\s A_2 A_1 A_1
    +\int_\s A_1 A_2 A_1+\cdots.
\]

The integral over $\s$ is the double integral of the pull-back to
$\Delta^2$. The factors of $A_1$ will just give the parallel
transport $\Phi_1$ along paths connecting $v_0$ and $v_2$ to the
point $v=\s(x,y)$
\[
    \psi_2(\s)=\int_{1\geq x\geq y\geq0} \s^\ast (\Phi_1(v_0,v)
    A_2(v)\Phi_1(v,v_2))
    \in \Hom(C(v_2),C(v_0))
\]
where $\Phi_1(v_0,v), \Phi_1(v,v_2)$ are given by parallel
transport along paths given by two straight lines each as shown in
the Figure.
\def\arrowheadup{\qbezier(0,.1)(0,.05)(.05,0)
      \qbezier(0,.1)(0,.05)(-.05,0)}
\def\arrowheaddown{\qbezier(0,0)(0,.05)(.05,0.1)
      \qbezier(0,0)(0,.05)(-.05,0.1)}
\def\arrowheadleft{\qbezier(0,0)(.05,0)(.1,.05)
      \qbezier(0,0)(.05,0)(.1,-.05)}
\def\arrowheadright{\qbezier(0.1,0)(.05,0)(0,.05)
      \qbezier(0.1,0)(.05,0)(0,-.05)}
\def\arrowheaddownleft{\qbezier(0,0)(0.05,.05)(.04,0.12)
      \qbezier(0,0)(0.05,.05)(.12,0.04)}
\begin{figure}[htbp]\label{fig1}
\caption{$\Phi_1(v_0,v), \Phi_1(v,v_2)$ are parallel transport
along dark lines.}
\begin{center}
%
\setlength{\unitlength}{1in}
\mbox{
\begin{picture}(4,1.5)
      \thinlines
      \put(1,0){      
      \put(1.25,.5){$v=\s(x,y)$}
	\put(0.3,0.1){$v_0$}
	\put(0.5,.2){\line(1,0){1}} 
	\put(1.5,.2){\line(0,1){1}}  
	\qbezier(0.5,.2)(1,.7)(1.5,1.2)  
	     \put(0.9,0){$\s(x,0)$}
	     \thicklines
	          \put(1.15,.85){\circle*{.05}}  
	               \put(1.15,.2){\line(0,1){.65}} 
	          \put(1.15,.2){\circle*{.05}}  
	               \put(1.15,.53){\circle*{.05}}
	           \put(.8,.2){\arrowheadleft}
      \put(.5,.2){\line(1,0){.65}}  
      \qbezier(1.15,.85)(1.3,1)(1.5,1.2) 
      \put(.5,.2){\circle*{.05}}  
      \put(1.5,1.2){\circle*{.05}}  
      \put(1.5,0.2){\circle*{.05}}  
           \put(.7,.9){$\s(x,x)$}
           \put(1.5,1.3){$v_2$}
           \put(1.6,0.1){$v_1$}
           \put(1.25,.95){\arrowheaddownleft}
      }
\end{picture}
}
\end{center}
\end{figure}
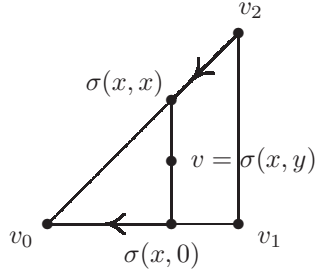
%
\[
    \Phi_1(v,v_2)=\Phi_1(v,\s(x,x))\Phi_1(\s(x,x),v_2)
\]
\[
    \Phi_1(v_0,v)=\Phi_1(v_0,\s(x,0))\Phi_1(\s(x,0),v).
\]

\begin{prop}\label{lemma 2}
$A_0(v_0)\psi_2(\s)+\psi_2(\s)A_0(v_2)=\Phi_1(v_0,v_2)-\Phi_1(v_0,v_1)\Phi_1(v_1,v_2)$.
\end{prop}

\begin{proof}
For $x\in[0,1]$ let $\Phi(x)$ be the parallel transport of
$d-A_1$ along the three segment path:
\[
    \Phi(x)=\Phi_1(v_0,\s(x,0))\Phi_1(\s(x,0),\s(x,x))\Phi_1(\s(x,x),v_2).
\]
Then
\[
    \Phi(0)=\Phi_1(v_0,v_2)
\]
\[
    \Phi(1)=\Phi_1(v_0,v_1)\Phi_1(v_1,v_2).
\]
So the right hand side of the formula we are proving is
\begin{equation}\label{eq:RHS of lemma 2}
    \Phi(0)-\Phi(1)=-\int_0^1 d\Phi(x).
\end{equation}
By Proposition \ref{lemma 1} the left hand side is equal to
\[
    \int_\s
    \Phi_1(v_0,v)[A_0(v)A_2(v)+A_2(v)A_0(v)]\Phi_1(v,v_2)
\]
(The sign in front of $A_2(v)A_0(v)$ is $(-1)^2=+1$ since the form degree of $A_2$ is 2.) By (\ref{eq:infinitesmal twisting cochain}) this is equal to
\begin{equation}\label{eq:LHS of lemma 2}
    =\int_\s
    \Phi_1(v_0,v)[-A_1(v)A_1(v)+dA_1(v)]\Phi_1(v,v_2)
\end{equation}

In (\ref{eq:RHS of lemma 2}), we have
\[
  -  \frac{d\Phi(x)}{dx}=\Phi_1(v_0,\s(x,0))X(x)\Phi_1(\s_1(x,x),v_2)
\]
where
\[
    X(x)=A_1^x\Phi_1(\s(x,0),\s(x,x))
    -\frac{d}{dx}\Phi_1(\s(x,0),\s(x,x))
    -\Phi_1(\s(x,0),\s(x,x))A_1^x
\]
using the notation $\s^\ast(A_1)=A_1^xdx+A_1^ydy$. (The term $\Phi_1(\s(x,0),\s(x,x))A_1^y$ which occurred with positive sign in the second term and negative sign in the third term was cancelled.) Comparing this to (\ref{eq:LHS of lemma 2}) we are reduced to showing that $X(x)=Y(x)$ where
\[
    Y(x)=\int_{0\leq y\leq
    z}dy\,\Phi_1(\s(x,0),v)\left(-A_1^xA_1^y+A_1^yA_1^x+\frac{\d A_1^y}{\d
    x}-\frac{\d A_1^x}{\d y}\right)\Phi_1(v,\s(x,x)).
\]

Expressing $\Phi_1(\s(x,0),\s(x,x))$ as an iterated integral of
$A_1^ydy$ we see that the second term of $X(x)$ is equal to the
third term of $Y(x)$ (with $\frac{\d A_1^y}{\d x}$). The negative sign comes from the fact that we are going backwards along the $y$ direction ($dt=-dy$). The other
three terms of $Y(x)$ form the commutator of $A_1^x$ with each
factor $A_1/dt=-A_1^y$ in the iterated integral representation of $\Phi_1(\s(x,0),\s(x,x))$. This can be more easily seen in the product limit form:
\[
	\Phi_1(\s(x,0),\s(x,x))=\lim_{\begin{matrix}n\to \infty\\\Delta y={x}/n\end{matrix}} \prod_{1}^{n} (I-A_1^y \Delta y)
\]
The commutator of $A_1^x$ with $I-A_1^y\Delta y$ is:
\[
	A_1^x(I-A_1^y\Delta y)-(I-A_1^y\Delta y)(A_1^x+\Delta_yA_1^x)=-A_1^xA_1^y\Delta y+A_1^y \Delta yA_1^x-\Delta_yA_1^x+o(\Delta y)
\]
So, the commutator of $A_1^x$ with $\Phi_1(\s(x,0),\s(x,x))$ is
\[
	A_1^x\Phi_1-\Phi_1A_1^x=\lim_{\Delta y\to 0} \sum_i \prod_{1}^{i-1} (I-A_1^y \Delta y)\left[-A_1^xA_1^y\Delta y+A_1^y \Delta yA_1^x-\Delta_yA_1^x\right]\prod_{i+1}^n (I-A_1^y \Delta y)
\]
So, the sum of the remaining two terms of $X(x)$ is equal to the sum of the remaining three terms of $Y(x)$ and we conclude that $X(x)=Y(x)$ proving the proposition.
\end{proof}

The construction that we just explained in detail is a special case of a construction outlined by Chen in \cite{Chen73}, sec. 4.5. Chen constructs mappings
\[
	\theta_{(n)}:I^{n-1}\to P(\Delta^n,v_n,v_0)
\]
from the $n-1$ cube $I^{n-1}$ to the space $P(\Delta^n,v_n,v_0)$ of smooth paths in $\Delta^n$ from $v_n$ to $v_0$ by smoothing a piecewise linear construction very similar to the one we explained. When $n=2$, this is the 1-parameter family of paths in $\Delta^2$ given in Figure 1 above. The main conceptual difference between Chen's construction and ours is that Chen follows this mapping with a smooth mapping of the $n$ simplex into the space $B$ which sends all vertices of $\Delta^n$ to the base point of $B$. He uses this to obtain a cubical chain complex for $\Omega B$.

A very longwinded description of the higher steps in this process of converting a flat superconnection into a simplicial twisting cochain can be found in \cite{IgIterated} which uses much of Chen's notation to allow for comparison and to make it easier to understand Chen's work.

In the last section of this paper we will show how the entire process can be done in an easier way.

%
%

\section{Another method}\label{sec8: another method}

There is another method for constructing a simplicial twisting
cochain from a flat connection. We assume that $B$ is compact and
we choose a finite ``good cover'' for $B$. (See \cite{BottTu}.)
This is a covering of $B$ by contractible open sets $U$ so that
all nonempty intersections
\[
    U_{\a_1}\cap U_{\a_2}\cap\cdots\cap U_{\a_n}
\]
are also contractible.

\begin{lem}\label{lem:superconnection complex on U}
If $U$ is a contractible open subset of $B$ and $D$ is a flat
superconnection on a graded vector bundle $V$ over $B$ then the
cohomology of the superconnection complex over $U$ is isomorphic
to the cohomology of $V$ using $A_0$ as differential:
\[
    H^n(\Omega(U,V|U),D)\cong H^n(V,A_0)
\] where the isomorphism is given by restriction to any point in
$U$.
\end{lem}

\begin{proof}
The spectral sequence collapses since its $E_2$-term is
$H^{p}(U;H^q(V))$.
\end{proof}

This lemma implies that
\[
    F:U\mapsto(\Omega(U,V|U),D)
\]
is a functor from the nerve $\cN\bu \cU$ of the good cover $\cU$ of $B$
to the category of cochain complexes over $\CC$ and cochain
homotopy equivalences. Applying the $A_\infty$ cohomology functor
we get a contravariant $A_\infty$ functor $H^\ast F$ on $\cN\bu\cU$. By Proposition \ref{dual of Ainfty is Ainfty} we can dualize to get a covariant $A_\infty$ functor $\Phi_1+\psi$ on $\cN\bu\cU$ with coefficients in $(\Phi,\Phi_1)=(V^\ast,\nabla^\ast)$. Subtracting $\Phi_1=\nabla^\ast$ we get the twisting cochain $\psi$ satisfying the following.

\begin{thm}\label{thm:superconnection complex is dual to twisted tensor product} The twisted tensor product $C_\ast(\cN\bu\cU)\otimes_\psi V^\ast$ is homotopy equivalent to the dual of the superconnection complex. I.e.,
\[
    (\Omega(B,V),D)\simeq \Hom(C_\ast(\cN\bu\cU)\otimes_\psi
    V^\ast,\CC).
\]
\end{thm}

\begin{proof}
This holds by induction on the number of open sets in the finite
good covering $\cU$. When the number is $1$ we use Lemma
\ref{lem:superconnection complex on U}. To increase the number we
use Mayer-Vietoris.
\end{proof}

\begin{rem} By Ed Brown's Theorem \ref{thm:Brown}
this implies that, if the superconnection is constructed
correctly, the superconnection complex gives the cohomology of the
total space of a smooth manifold bundle. This construction also
allows us to compare flat superconnections with twisting cochains,
giving a K-theory difference class with a well defined higher
torsion. But, this is the subject of another paper.
\end{rem}


\end{document}